\documentclass[12pt]{amsart}
\usepackage{amsmath}
\usepackage{amsthm}
\usepackage{mathrsfs}
\usepackage{amsfonts} 
\usepackage{graphicx}
\usepackage{hyperref}
\usepackage{epsfig}

\numberwithin{equation}{section}

\allowdisplaybreaks

\newtheorem{theorem}{Theorem}[section]

\newtheorem{proposition}[theorem]{Proposition}
\newtheorem{lemma}[theorem]{Lemma}
\newtheorem{remark}[theorem]{Remark}
\newtheorem{definition}[theorem]{Definition}
\newtheorem{corollary}[theorem]{Corollary}

\newcommand{\ep}{\varepsilon}

\newcommand{\R}{\mathbb{R}}
\newcommand{\N}{\mathbb{N}}
\newcommand{\Z}{\mathbb{Z}}
 
\newcommand{\PR}{\mathbb{P}}

\newcommand{\dist}{\vert\vert}

\newcommand{\visit}{\mathcal{V}}

\newcommand{\sab}{\mathsf{P_{SAB}}}

\newcommand{\isabstar}{\mathsf{P^*_{iSAB}}}
\newcommand{\saw}{\mathsf{P_{SAW}}}

\newcommand{\stickbreak}{\mathsf{StickBreak}}

\newcommand{\dr}{\mathsf{D}}

\newcommand{\sawset}{{\rm SAW}}
\newcommand{\sabset}{{\rm SAB}}

\newcommand{\isabset}{{\rm iSAB}}

\newcommand{\zigzag}{\mathsf{ZigZag}}

\newcommand{\reflect}{\mathcal{R}}

\newcommand{\unfold}[3]{{{\mathsf{Unf}_{#1,#2}(#3)}}}
\newcommand{\unf}{\mathsf{Unf}}
\newcommand{\multiunf}{\mathsf{MultiUnf}}

\newcommand{\latt}{\Z^d}
\newcommand{\lattv}{\Z^d}
\newcommand{\latte}{E(\Z^d)}

\newcommand{\muc}{\mu_c}

\newcommand{\shortzz}{\mathsf{ShortZZ}}

\newcommand{\mzfr}{\mathrm{ManyZFewR}}

\newcommand{\ner}{\mathsf{R}}

\newcommand{\prne}{\PR_{{\rm iSAB}}}
\newcommand{\prneinf}{\PR_{{\rm iSAB}}^{\otimes\mathbb N}}
\newcommand{\prneinfbi}{\PR_{{\rm iSAB}}^{\otimes \mathbb Z}}

\newcommand{\calO}{\mathcal O}

\newcommand{\an}{\widehat{an}}
\newcommand{\bn}{\widehat{bn}}

\title{Self-avoiding walk is sub-ballistic
}
\date{}

\author[H.~Duminil-Copin and A.~Hammond]{Hugo Duminil-Copin and Alan Hammond}
\address{ Department of Statistics, 
           University of Oxford,
                  1 South Parks Road,
                 Oxford, OX1 3TG, U.K.}
                 
                 \address{D\'epartement de Math\'ematiques, Universit\'e de Gen\`eve, 2--4 rue du Li\`evre, Gen\`eve, Switzerland.
} 
\email{duminil.copin@unige.ch,hammond@stats.ox.ac.uk
}

\begin{document}

\begin{abstract}
We prove that self-avoiding walk on $\Z^d$ is sub-ballistic in any dimension $d \geq 2$. That is, writing $\dist u \dist$ for the Euclidean norm of $u \in \Z^d$, and $\saw_n$ for the uniform measure on self-avoiding walks $\gamma:\{0,\ldots,n\} \to \Z^d$ for which $\gamma_0 = 0$, we show that, for each $v > 0$, there exists $\ep > 0$ such that, for each $n \in \N$, $\saw_n \big(  \max\big\{\dist \gamma_k \dist : 0 \le k \le n\big\} \geq v n \big) \leq e^{-\ep n}$.
\end{abstract}

\maketitle

\section{Introduction}

Flory and Orr \cite{Flory,Orr47} introduced self-avoiding walk as a simple model of a polymer. It quickly became clear that the model exhibits rich behaviour the rigorous understanding of which poses a real challenge to mathematicians.
Originally, Flory was interested in the typical behaviour of self-avoiding walk, and in particular in the mean squared displacement of its endpoint. Since then, much effort has been invested in the study of typical geometric properties of the model.

\subsection{The model and the results} Let $d \geq 2$. For $u \in \R^d$, let $\dist u \dist$ denote the Euclidean norm of $u$. 
Let $\latte$ denote the set of nearest-neighbour bonds of the integer lattice $\latt$.
A {\em walk} of length $n \in \N$ is a map $\gamma:\{0,\ldots,n \} \to \lattv$ such that $(\gamma_i,\gamma_{i+1}) \in \latte$ for each $i \in \{0,\ldots,n-1\}$.
An injective walk is called {\em self-avoiding}. 
Let $\sawset_n$ denote the set of self-avoiding walks of length~$n$ that start at $0$ and let $\saw_n$ denote the uniform law on $\sawset_n$. 
This article is devoted to proving that self-avoiding walk is sub-ballistic:
\begin{theorem}\label{sub SAW}
Let $v>0$. There exists $\ep> 0$ such that, for each $n \in \N$,
$$
 \saw_n \Big( \max \big\{\dist \gamma_k \dist : 0 \le k \le n \big\} \ge v n \Big) \le e^{- \ep n} \, .
$$
\end{theorem}

The theorem has the following immediate consequence.
\begin{corollary}\label{cor:1}
Set  $\langle \dist \gamma_n \dist^2 \rangle = \tfrac{1}{|\sawset_n|}\sum_{\gamma\in \sawset_n} \dist \gamma_n \dist^2$. Then
$$
\lim_{n\rightarrow \infty} n^{-2} \langle \dist \gamma_n \dist^2 \rangle = 0 \, .
$$
\end{corollary}

\subsection{Conjectures on the mean-square displacement}
Numerical computations and non-rigorous theory predict the precise behaviour of the mean-squared displacement of the walk's endpoint. It is conjectured that
$$\langle \dist \gamma_n \dist^2 \rangle~=~n^{2\nu+o(1)}\text{ where }\nu=\begin{cases}1 & d=1\\
3/4& d=2\\
\approx 0.59\cdots & d=3 \\
1/2 &d=4\\
1/2 & d\ge 5 \, . \end{cases}$$

The behaviour predicted by this conjecture should be compared to that of simple random walk. 
In general, the displacement of the endpoint of self-avoiding walk is expected to exceed that of simple random walk.  
In dimensions two and three, this difference is manifested in a strong form, with the value of $\nu$ in the self-avoiding case exceeding its counterpart for the simple one (which is $1/2$). 
Dimension four is known as the {\em upper critical dimension}: here, the two values of $\nu$ coincide, with self-avoiding walk experiencing further displacement in the form of a poly-logarithmic correction. In dimensions five and higher, the two values coincide and no logarithmic corrections occur; indeed, the scaling limit of each process is Brownian motion. However, the diffusion rates of the Brownian motions in the simple random walk and self-avoiding walk scaling limits differ, with the rate being higher in the self-avoiding case.

In dimensions five and above, the conjecture was proved for a version of self-avoiding walk with weak repulsion in
Brydges and Spencer \cite{BS85}. Hara and Slade
\cite{HS91,HS92} proved the conjecture for self-avoiding walk and established that this process has Brownian motion as a scaling limit in these dimensions.
% (see also the book \cite{MS93}). 

Rigorous analysis in dimension four is much more subtle.
Recently some impressive results have been achieved using
a supersymmetric renormalization group approach. These results concern continuous-time weakly self-avoiding walk: see
\cite{BIS09,BS10,BDGS11} and references within. 

As far as we know, no significant physical predictions have been made regarding the scaling limit of self-avoiding walk in three dimensions.
In dimension two, the Coulomb gas formalism  \cite{Nie82,Nie84} provides the prediction that $\nu = 3/4$, with this prediction later being made \cite{Dup89,Dup90} using conformal field theory. The scaling limit of the model in this dimension has been identified~\cite{LSW5} subject to certain assumptions principal among which is conformal invariance: given that the scaling limit may be expected to be supported on simple curves and to enjoy a natural restriction property, the assumption of a certain conformal covariance forces the limit to be given by Schramm-Loewner Evolution with parameter $\kappa = 8/3$. 

There are however very few rigorous unconditional results regarding dimensions two and three. The question of a non-trivial upper bound on mean-squared endpoint displacement is raised in the introduction of \cite{MS93}. That the intuitively very natural assertion that self-avoiding walk in dimensions two and three is sub-ballistic has remained unresolved for a long time is one of several archetypical examples which bear witness to  the theoretical difficulty that this model presents in the low-dimensional case.

\subsection{Self-avoiding walk between two points of a domain}
\begin{figure}
\begin{center}
\includegraphics[width=0.50\textwidth]{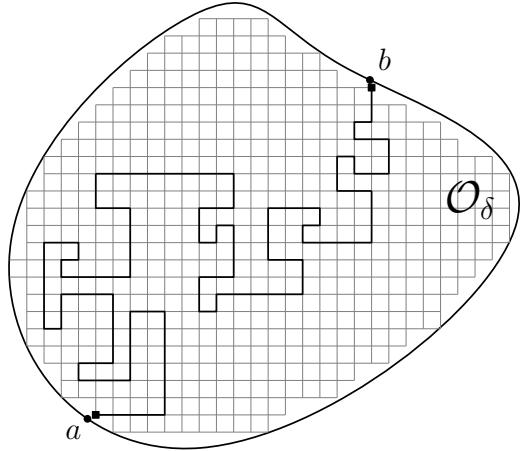}
\end{center}
\caption{\label{fig:D_delta}
The path $\gamma_\delta$ in $(\mathcal O_\delta,a_\delta,b_\delta)$.}
\end{figure}

The self-avoiding walk model exhibits a phase transition when defined slightly differently. Let $\calO$ be a simply connected smooth domain in $\mathbb R^d$ with two points $a,b$ on the boundary. For $\delta>0$, let $\calO_\delta$ be the largest
connected component of $\calO \cap \delta\mathbb Z^d$ and let $a_\delta$, $b_\delta$ be the two sites of $\calO_\delta$ closest to $a$ and $b$ respectively. We think of $(\calO_\delta,a_\delta,b_\delta)$ as being an approximation of
$(\calO,a,b)$; see Figure~\ref{fig:D_delta}. Let $z>0$. On $(\calO_\delta,a_\delta,b_\delta)$, define a probability measure on the finite set of self-avoiding walks in $\calO_\delta$ from $a_\delta$ to $b_\delta$ by the formula
\begin{equation} 
\mathbb P_{(\calO_\delta,a_\delta,b_\delta,z)}(\gamma_\delta)~=~\frac{z^{|\gamma_\delta|}}{Z_{(\calO_\delta,a_\delta,b_\delta)}(z)} \, ,
\end{equation}
where $|\gamma_\delta|$ is the length of $\gamma_\delta$, and $Z_{(\calO_\delta,a_\delta,b_\delta)}(z)$ is a normalizing factor.

%\begin{figure}
%\centering
%\includegraphics[width= 6cm,bb= 69 192 339 460,clip]{Ddelta.eps}
%\caption{\label{fig:D_delta}
%The path $\gamma_\delta$ in $(\mathcal O_\delta,a_\delta,b_\delta)$.}
%\end{figure}

A phase transition occurs at the value $\mu_c^{-1}$, where $\mu_c$ is the connective constant (whose definition we will shortly provide). For $z<\mu_c^{-1}$, there exists $C=C(z)>0$ such that $\mathbb P_{(\calO_\delta,a_\delta,b_\delta,z)}(|\gamma_\delta|\le C/\delta)\rightarrow 1$ as $\delta \searrow 0$. In fact, $\gamma_\delta$ has Gaussian fluctuation of
order $\sqrt \delta$ around the geodesic between $a$ and $b$
in $\calO$ (assuming the geodesic is unique --- otherwise some changes are needed). Broadly, these results may be expected to be a consequence of the theory developed by Ioffe \cite{Ioffe} that treats unrestricted self-avoiding
walk in $\Z^d$; however, they have not been rigorously derived to the best of our knowledge. For $z>\mu_c^{-1}$, the converse is true \cite{DKY11} in the sense that $\gamma_\delta$ becomes space-filling.

Viewed in this light, Theorem~\ref{sub SAW} rules out the possibility that the model's critical behaviour coincides with that of the subcritical phase $z < \mu_c^{-1}$:
\begin{corollary}\label{cor:2}
Let $(\calO,a,b)$ be such that $\partial \calO$ is smooth in a neighbourhood of $a$ and of $b$. For every $K>0$, $\mathbb P_{(\calO_\delta,a_\delta,b_\delta,\mu_c^{-1})}\big(|\gamma_\delta|\le K/\delta\big) \to 0$ as $\delta \searrow 0$.
\end{corollary}

%In \cite{LSW5}, the exponent in dimension two was computed using the strong hypothesis that when $x=x_c$, the curve $\gamma_\delta$ converges to a random continuous simple curve, called the Schramm-L\"owner Evolution of parameter $8/3$; see \cite{LSW5} and references therein for more information on the Schramm-L\"owner Evolution.

\subsection{More general graphs}
The proof of Theorem~\ref{sub SAW} can be extended (with additional technicalities) to lattices with symmetry. For instance, the case of the hexagonal lattice should follow from the same reasoning. While the question can be asked on any locally-finite infinite transitive graph, the answer will differ drastically depending on the growth of the graph. For instance, it is strongly expected that self-avoiding walk on the Cayley graph of a non-amenable group is ballistic (see the upcoming Question $5$); see \cite{MW05} for a study of some hyperbolic planar graphs.

In \cite{BCKM00} is considered a model where a walk $\gamma \in \sawset_n$ is chosen with weight proportional to $\prod_{i=1}^d z_{i,+}^{N_{i,+}} z_{i,-}^{N_{i,-}}$ where $N_{i,+/-}$ is the number of increasing or decreasing steps made by $\gamma$ in direction $e_i$ and $z_{i,+/-} \in (0,\infty)$ are parameters. In the asymmetric case that $z_{i,+} \not= z_{i,-}$ for some $i \in [1,d]$ (and no $z_{j,+/-} = 0$), \cite[Theorem 1.2]{BCKM00} is a result to the effect that the walk is typically ballistic. 
\medbreak
\subsection{Open problems}
\medbreak
Theorem~\ref{sub SAW} is a first step towards proving the conjecture on mean-squared displacement. An interesting improvement would be a quantitative bound: \medbreak 
\noindent\textbf{Question 1.} Show that for some $\ep>0$, $\lim_{n\rightarrow \infty} n^{-(2-\ep)} \langle \dist \gamma_n \dist^2 \rangle =0$.
\medbreak
Theorem~\ref{sub SAW} raises the prospect of confirming that the walk from the lower-left to upper-right corner of a square sways to macroscopic distance from the diagonal, but it does not resolve this question: 
\medbreak
\noindent\textbf{Question 2.} Fix $(\calO,a,b)$. Show that $\gamma_\delta$ does not converge to a geodesic.
\medbreak
The theorem rules out the extreme of fast movement by the walk. What about the other extreme? Does self-avoiding walk move further than simple random walk?
\medbreak
\noindent\textbf{Question 3.} Show that $\liminf_{n\rightarrow \infty} n^{-1} \langle \dist \gamma_n \dist^2 \rangle > 0$. 
\medbreak
In fact, far slower motion by self-avoiding walk has yet to be ruled out. In principle, the walk may be space-filling: 
\medbreak
\noindent\textbf{Question 4.} Show that $\lim_{n\rightarrow \infty} n^{-2/d} \langle \dist \gamma_n \dist^2 \rangle =\infty$. 
\medbreak
As we have mentioned, the question of ballisticity may be posed for self-avoiding walk on many transitive graphs, such as Cayley graphs. Let $G$ be a finitely generated infinite group and let $S$ be a finite symmetric system of generators. Let $\mathcal G_{G,S}$ be the Cayley graph associated to $(G,S)$. 
 \medbreak
\noindent\textbf{Question 5.} Is self-avoiding walk ballistic whenever simple random walk is? In particular, is the model ballistic if $G$ is non-amenable?

\subsection{Acknowledgments.} The authors would like to thank Tony Guttmann, Neal Madras and Yvan Velenik for very careful readings of the article, and for many valuable comments. H.D-C. was supported by ANR grant BLAN06-3-134462,
the EU Marie-Curie RTN CODY, the ERC AG CONFRA, as well as by the Swiss
{FNS}. A.H. was supported by EPSRC grant EP/I004378/1 and thanks the Mathematics Department of the University of Geneva for its hospitality during the year 2011-12.

\section{Preliminaries}

\subsection{Notation}
The symbol $\mathbb N$ denotes the set of integers $\{0,1,2,3,\ldots\}$. We set $[a,b]=\{a,a+1,\dots,b\}$. 

For $u \in \R^d$, we write $u = (u_1,\ldots,u_d)$. We also write $x(u) = u_1$ and $y(u) = u_2$.
For $k \in [1,d]$,
let $e_k$ be the vector of Euclidean norm one whose $k^{\textrm{th}}$ entry is $1$.  As we tend to visualise our constructions for the two-dimensional case, we will sometimes refer to the directions $e_1$ and $e_2$ as east and north.

The cardinality of a finite set $A$ is denoted by $|A|$. Abusing slightly this notation, the length of a walk $\gamma$ will be denoted by $|\gamma|$; recall that $| \gamma |$ is the number of edges (rather than the number of vertices) comprising $\gamma$.
\medbreak
For $m,n \in \N$, let $\gamma$ and $\tilde{\gamma}$ be two walks of lengths $m$ and $n$, neither of which need to start at $0$. The concatenation $\gamma \circ \tilde\gamma$ of $\gamma$ and $\tilde \gamma$ is given by 
$$
\big( \gamma \circ \tilde\gamma \big)_k = \begin{cases}
\gamma_k & k\le m \, ,\\
\gamma_m \, + \,  \big(\tilde \gamma_{k-m}-\tilde\gamma_0 \big) & m+1\le k\le m+n \, .
\end{cases}
$$
\medbreak
Let $A$ and $B$ be two sets. Let $\mathfrak{P}(B)$ denote the power-set of a set $B$. A multi-valued map from $A$ to $B$ is a function $f:A\rightarrow \mathfrak{P}(B)$. 
We will use such maps in order to estimate the size of sets. On several occasions, we will use the following basic principle, which we name for future reference:
\medbreak
\noindent\textbf{Multi-valued map principle.} {\em Let $f$ be a multi-valued map from $A$ to $B$. Write $f^{-1}(b):=\{a\in A:b\in f(a)\}$. Then
$$
 |A| \le \frac{\max_{b\in B}|f^{-1}(b)|}{\min_{a\in A}|f(a)|} \cdot |B| \, .
$$}
\subsection{Bridges}
Recall the classical definition of a bridge. A self-avoiding walk is called a {\em self-avoiding bridge} if
\begin{itemize}
 \item 
its first element attains uniquely the minimal  $y$-coordinate on the walk;
 \item and
its final element attains, not necessarily uniquely, the maximal  $y$-coordinate on the walk.
\end{itemize} 
The self-avoiding walk being the object of attention, we will usually omit the term ``self-avoiding'' in referring to walks and bridges.
For $n \in \N$, let $\sabset_n$ and                                                                                  $\sawset_n$ denote the set of bridges and walks of length $n$. Let $\sawset$ and $\sabset$ be the corresponding sets where now no condition on the length is applied.
Let $\sab_n$ and $\saw_n$ denote the uniform law on $\sabset_n$ and $\sawset_n$. 

\begin{figure}[ht]
\centering
\epsfig{file=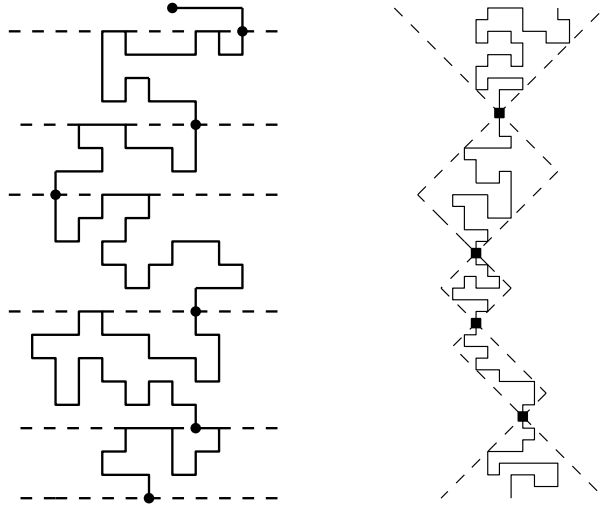, width=8cm}
\caption{\label{fig:ramp}
On the left, a bridge with its renewal points. On the right, a bridge with its diamond points (which we will shortly discuss).}\label{ramprenewal}
\end{figure}

We begin by recording a straightforward statement concerning the number of bridges of length $n$. Recall from \cite{hammersleywelsh} that a standard subadditive argument assures the existence of the connective constant $\mu_c$, given by $\mu_c =  \lim_n \big| \sawset_n \big|^{1/n} \in (0,\infty)$. (See \cite{DS10} for a rigorous identification of $\muc$'s value for the hexagonal lattice.)

\begin{proposition}\label{number ramps} 
The generating functions of walks and bridges satisfy
\begin{align}
\sum_{\gamma\in\sawset}z^{|\gamma|}&\le z^{-1}\exp \Big(2\sum_{\gamma\in\sabset}z^{|\gamma|}-2\Big) \, .
\label{eq:generating bridge}
\end{align}
Furthermore, there exists $C>0$ such that, for each $n \in \N$,
\begin{equation}\label{eq:22}
e^{-C\sqrt n} \vert \sawset_n \vert \le |\sabset_n|\le  \vert \sawset_n \vert \, .
\end{equation}
\end{proposition}

In a classic argument of Hammersley and Welsh \cite{hammersleywelsh}, (\ref{eq:22}) is demonstrated by the construction of a map from 
 $\sawset_n$ and $\sabset_n$ enjoying the property that each preimage has cardinality at most $e^{C \sqrt{n}}$; formally, then, the multi-valued map principle (applied in this case to single-valued maps) yields~\eqref{eq:22}. In fact, the techniques imply \eqref{eq:generating bridge}: see \cite[(3.1.13)]{MS93}.
\medbreak
The notion of cutting an object into irreducible pieces has been much used, largely because it is crucial for developing a renewal theory.  This theory was developed for bridges in \cite{kestenone,kestentwo}.
The set $\ner_\gamma$ of {\em  renewal points}  of $\gamma\in \sabset_n$ is the set of points of the form $\gamma_i$ with $i\in[0, n]$, for which $\gamma[0,i]$ and $\gamma[i,n]$ are  bridges.  %When no confusion is possible, we write $\ner $ for $ \ner_\gamma$. 
See Figure~\ref{fig:ramp}.
A  bridge $\gamma\in \sabset_n$ with $n\ge 1$ is said to be {\em irreducible} if $\gamma_k$ is not a renewal point for any $k\in [1,n-1]$.

Let $\isabset$ be the set of irreducible  bridges of arbitrary length whose initial point is $0$. Every  bridge is the concatenation of a finite number of irreducible elements; the decomposition is unique and the set $\ner_\gamma$ is the union of the initial and terminal points of the  bridges that comprise this decomposition. This leads to the following result due to Kesten.

\begin{lemma}[Kesten] \label{Kesten relation}
We have that  
$\displaystyle
\sum_{\gamma\in \isabset}\mu_c^{-|\gamma|}=1.$%\quad\text{ and }\quad\sum_{\gamma\in \isar}\mu_c^{-|\gamma|}=1 \, .$$
\end{lemma}
\noindent{\bf Proof.}
Since the argument is short, we provide it here. The decomposition of walks into irreducible  bridges yields
$$\sum_{\gamma\in\sabset}z^{|\gamma|}=\frac{1}{1-\sum_{\gamma\in \isabset}z^{|\gamma|}}$$
for every $z$ such that the two series converge. Proposition~\ref{number ramps} shows that the radius of convergence of the series on the left-hand side  is $\mu_c^{-1}$. The generating function does not blow up for $z<\mu_c^{-1}$, and thus $\sum_{\gamma\in \isabset}z^{|\gamma|}<1$ for any $z<\mu_c^{-1}$, whence $\sum_{\gamma\in \isabset}\mu_c^{-|\gamma|}\le 1$ follows by means of the monotone convergence theorem. 

Furthermore, this sum will equal one if $\sum_{\gamma\in\sabset}z^{|\gamma|}$ diverges as $z$ tends to $\mu_c^{-1}$. But \eqref{eq:generating bridge} shows that $\sum_{\gamma\in\sabset}z^{|\gamma|}$ diverges as $z\nearrow \mu_c^{-1}$ if $\sum_{\gamma\in{\rm SAW}}z^{|\gamma|}$ diverges in this limit. Divergence of the latter sum follows from  $\vert \sawset_n \vert \geq \muc^n$, which is a consequence of the submultiplicative bound $\vert \sawset_{n+m} \vert \leq \vert \sawset_n \vert\vert \sawset_m \vert$. \qed

By means of Lemma~\ref{Kesten relation}, we define the probability measure $\prne$ on $\mathrm{iSAB}$ by setting 
$\prne(\gamma) = \mu_c^{-\vert \gamma \vert}$. 

\subsection{Statements and sketches of the proofs of the main elements}\label{secsketch}

We will prove Theorem \ref{sub SAW} by contradiction. We might argue instead by conditioning but the further notation required would be substantial.

We first work with  bridges. We record the assumption whose negation we seek to prove by contradiction; when a result requires the assumption, we will say so in its statement.
\medbreak
\noindent{\bf Ballistic assumption.} {\em Suppose that for some $v > 0$, 
$$
\limsup_{n\rightarrow \infty} n^{-1}   \log \sab_{n} \Big(   y(\gamma_n)   \geq v n  \Big) =0 \, .
$$}
We first prove that, under the ballistic assumption,  the probability that a bridge possesses a positive density of renewal points decays subexponentially. More precisely: 
\begin{theorem}\label{thm pos dens}
Suppose that the ballistic assumption holds. Then there exists $\delta>0$ such that
$$
\limsup_{n\rightarrow \infty} n^{-1} \log \sab_n \big( |\ner_\gamma|\ge \delta n \big) =0 \, .
$$
\end{theorem}

While its statement is very intuitive, the proof of this result contains some of the central ideas of this paper.  The fact that self-avoiding walk is not Markovian makes the study delicate. The argument begins by noting that, if a  bridge travels ballistically in the northerly direction, many sets of the form $\{u\in\mathbb Z^d:y(u)=h\}$ are visited on at most a few occasions. We then use unfoldings in order to prove that, by paying a subexponential cost,
 many of these sets are in fact visited only once.

We will see that the macroscopic presence of such  renewal points at subexponential cost implies the stronger statement that  bridges typically have a positive density of renewal points, or, in other words, that the average size of an irreducible  bridge is finite.
\begin{corollary}\label{crucial}
Supposing the ballistic assumption, $\mathbb E_{\rm iSAB}(|\gamma|) <\infty$.
\end{corollary}
An independent argument rules out this possibility.
\begin{theorem}\label{no renewal}
We have that $\mathbb E_{\rm iSAB}(|\gamma|) =\infty$.
\end{theorem}

We now explain roughly how we will prove Theorem~\ref{no renewal}. Seeking a contradiction, we assume that $\mathbb E_{\rm iSAB}(|\gamma|) < \infty$.

Let $\prneinf$ denote the law on semi-infinite walks $\gamma:\N \to \lattv$ formed by the concatenation of infinitely many independent samples of $\prne$. This measure may be viewed heuristically as the law of a semi-infinite self-avoiding walk conditioned to remain in the upper-half plane. 
See \cite[Section 8.3]{MS93} for further discussion about $\prneinf$.

Our assumption that the irreducible renewal block has finite mean length means that we may invoke the law of large numbers to conclude that, under the law $\prneinf$, the
walk $\gamma$ proceeds northwards at a constant rate $\nu \in (0,1)$, so that $y(\gamma_n) = \nu n \big( 1 + o_{n\rightarrow \infty}(1) \big)$ as $n \to \infty$.
 However, for a typical sample $\gamma$ of $\prneinf$, the interval $[0,\nu n]$ is populated to positive density by renewal points; 
there are thus an order of $n^2$ pairs of renewal points such that each element of the pair has $y$-coordinate in this interval and where the distance between the two points is at least a small constant multiple of $n$. Writing $(r,r')$ for the indices of such a pair of renewal points,
we may perform a surgery on $\gamma$, rotating $\gamma_{[r,\infty)}$ clockwise about $\gamma_r$ by a right-angle, advancing along the new curve   
 to the image of $\gamma_{r'}$, and then rotating the subsequent trajectory counterclockwise by a right-angle about the encountered point, so that the altered curve again points in the northerly direction. 
The new curve thus proceeds northwards, then turns eastward for at least small positive multiple of $n$ steps, before returning northward. 
To a typical sample of $\prneinf$ we thus associate an order of $n^2$ alternative walks. To each alternative walk, there are at most order $n^2$ candidates for the original walk to which surgery was applied: indeed, to perform this reconstruction, it is enough to find each of the renewal point indices $r$ and $r'$, and  a factor of at most $n$ is needed for each one.

Thus, by an argument in the spirit of the multi-valued map principle, were all such alternative walks self-avoiding, they would have a $\prneinf$-probability which is bounded away from zero uniformly in $n$. However, alternative walks reach $y$-coordinate $\nu n$ after a number of steps which exceeds a multiple of $n$ where the factor is a constant which strictly exceeds one; this is in contradiction to $y(\gamma_n)=n(1+o_{n\rightarrow \infty}(1))$ ensured by the assumption that $\mathbb E_{\rm iSAB}(|\gamma|) < \infty$. This is an outline of an argument to arrive at a contradiction which would prove Theorem~\ref{no renewal}.

\begin{figure}
\centering
\epsfig{file=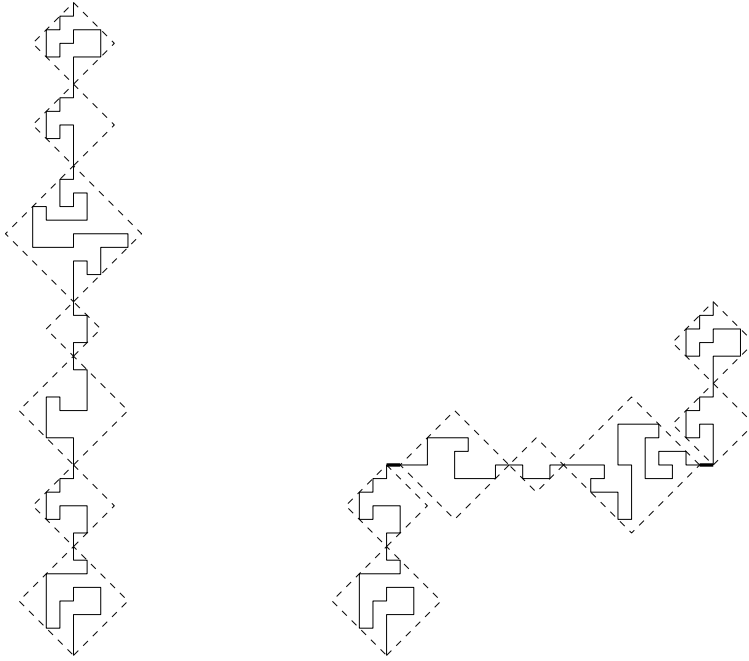, width=10cm}
\caption{The stickbreak operation is depicted. Note the two east-west edges, depicted in bold, that are inserted at either end of the rotated section in order to ensure that the outcome of the operation is self-avoiding.}\label{stickbreak}
\end{figure}

The difficulty with the plan is that our notion of renewal point is not sufficient to ensure that right-angle rotation about a renewal point of the path subsequent to that point results in a self-avoiding walk. Thus we introduce a stronger notion of renewal point for which this property is (in essence) valid. The point $\gamma_i$ is a {\em diamond point} of $\gamma$ if 
 \begin{itemize}
 \item for any $j\ge i$, $(x+y)(\gamma_j)\ge (x+y)(\gamma_i)$ and $(y-x)(\gamma_j)\ge (y-x)(\gamma_i)$,
 \item for any $j\le i$, $(x+y)(\gamma_j)\le (x+y)(\gamma_i)$ and $(y-x)(\gamma_j)\le (y-x)(\gamma_i)$.
 \end{itemize}
 Note that a diamond point is indeed a renewal point.

To implement our plan using diamond points in place of renewal points, we will first prove that, if $\mathbb E_{\rm iSAB}(|\gamma|) < \infty$, then 
a positive proportion of renewal points are in fact diamond points, so that the latter points also populate the $y$-coordinate to positive density.
Our operation of turning the curve first clockwise and then counterclockwise again at appropriate diamond points will be called the stickbreak operation and is illustrated in Figure~\ref{stickbreak}. Stickbreak now produces self-avoiding outputs and this allows us to derive the sought contradiction.

\subsection{Deducing Theorem~\ref{sub SAW}.}
Here, we apply %Theorem~\ref{thm pos dens}, 
Corollary~\ref{crucial} and Theorem~\ref{no renewal} to obtain Theorem~\ref{sub SAW}. 
\medbreak
\noindent{\bf Proof of Theorem~\ref{sub SAW}.} By the axial symmetry of the law $\saw_n$, note that  
\begin{eqnarray*}
 & & \saw_n \Big( \max\big\{\dist \gamma_k \dist: 0 \le k \le n \big\} \ge v n \Big) \\
 & \le & 2d \saw_n \Big( \max \big\{y(\gamma_k): 0 \le k \le n \big\}  \ge   d^{-1/2}v n \Big) \, .
\end{eqnarray*}

We apply the classical unfolding operation from walks into bridges (see \cite[Section 3.1]{MS93}) to walks $\gamma \in \sawset_n$ verifying  $\max\big\{y(\gamma_k): 0 \le k \le n \big\}  \ge   d^{-1/2}v n$.
Noting that $\max\big\{y(\gamma_k): 0 \le k \le n \big\}$ increases at each step, and that the outcome of the unfolding, being an element of $\sabset_n$, has the property that this maximum is attained by $k=n$, we learn that 
$$
\sab_n \Big( y(\gamma_n)   \geq d^{-1/2} v n \Big)\ge \tfrac{1}{2d}e^{-C\sqrt n}\saw_n \Big( \dist \gamma_n \dist \ge v n \Big) \, .
$$
Therefore, if the conclusion of Theorem~\ref{sub SAW} is violated, the ballistic assumption must be verified. Corollary~\ref{crucial} then contradicts Theorem~\ref{no renewal}.
\qed

\subsection{Structure of the paper}

Section~\ref{sec:ne} contains the proofs of Theorem~\ref{thm pos dens} and Corollary~\ref{crucial}.
In Section~\ref{sec:ir}, Theorem~\ref{no renewal} is proved. The final Section~\ref{seccoroll} is devoted to the proof of the corollaries.

\section{Positive density of  renewal points \\ at subexponential cost}\label{sec:ne}
In this section, we prove Theorem~\ref{thm pos dens}. For $v > 0$, set
\begin{equation}\label{eqsanenv}
\sabset_{n,v} = \Big\{ \gamma \in \sabset_n: y(\gamma_n) \geq vn  \Big\} \, . 
\end{equation}
Theorem~\ref{thm pos dens} follows from the ballistic assumption and the following result.
\begin{theorem}\label{thm pos dens 2}
Let $v > 0$. Let $\big\{ u_n: n \in \N \big\}$ be a subsequence of $\N$. 
There exists $\delta>0$ as well as 
a subsequence $\big\{ v_n: n \in \N \big\}$  of $\big\{ u_n: n \in \N \big\}$ and a sequence $\big\{ \ep_n: n \in \N \big\}$ with $\ep_n \searrow 0$ such that
$$
 \sab_{v_n,v} \Big(    \vert \ner_\gamma \vert    \geq \delta v_n \Big) \geq  e^{-\ep_n v_n} \, .
$$
\end{theorem}

For $h \in \Z$, let $E_h$ denote the set of (necessarily north-south) edges $e = (u,u') \in \latte$ such that $\big\{ y(u),y(u')\big\}  = \big\{ h, h + 1\big\}$.
For $\gamma \in \sawset_n$, define the {\it $h$-visiting edge-set}
 $\visit_{h,h+1}$ by $\visit_{h,h+1}(\gamma) =  E_h \cap \{(\gamma_i,\gamma_{i+1}): 0 \leq i \leq n-1 \big\}$.
%
%\begin{figure}
%\centering
%\epsfig{file=E_k.eps, width=5cm}
%\caption{\label{fig:E_k}The sets $E_h$ and $\visit_{h,h+1}$.}
%\end{figure}

 Let $\delta > 0$. For $m\ge 1$, write $\sabset_{n,v,\delta}^m \subseteq \sabset_{n,v}$ for the set of $\gamma \in \sabset_{n,v}$ such that there are at least $\delta n$ values of $h \in \Z$
for which $\vert \visit_{h,h+1} \vert \leq m$. Note that $\sabset_{n,v,\delta}^1 \subseteq \big\{\gamma\in\sabset_n: \vert \ner_\gamma \vert \geq \delta n  \big\}$ since if a walk $\gamma$ crosses from $E_h$ to $E_{h+1}$ only at $(\gamma_i,\gamma_{i+1})$, then $\gamma_i$ is a  renewal point. Thus, to derive Theorem~\ref{thm pos dens 2}, we may argue that $\sabset_{n,v,\delta}^1$ is not too small with respect to $\sabset_{n,v}$. 
The next proposition is the main technical result used in this deduction. 
It shows that, judged on an exponential scale and along a suitable subsequence, 
the sets $\sabset_{n,v,\delta}^k$ and $\sabset_{n,v,\delta'}^{k-1}$ have the same size if $\delta'\ll \delta$.

\begin{proposition}\label{pos dens}
Let $v > 0$. Let $k \geq 2$. Let $\delta>0$ and $\big\{ u_n: n \in \N \big\}$ be a subsequence of $\N$. 
Then there exist $\delta'>0$  and a subsequence  $\big\{ v_n: n \in \N \big\}$ of  $\big\{ u_n: n \in \N \big\}$ such that 
$$
\lim_{n\rightarrow \infty}\frac{1}{v_n}\log\bigg( \frac{|\sabset_{v_n,v,\delta}^k|}{|\sabset_{v_n,v,\delta'}^{k-1}|} \bigg) \, = \, 0 \, .
$$
\end{proposition}
\noindent{\bf Proof of Theorem \ref{thm pos dens 2}.} We claim that
\begin{equation}\label{sanenv inc}
\sabset_{n,v} \subseteq \sabset_{n,v,v/2}^{\lceil 2/v\rceil} \, .
\end{equation} 
To verify this,
note that $\gamma \in \sabset_{n,v}$ implies that $\visit_{h,h+1} \not= \emptyset$ for at least $v n$ values of $h \in \Z$. Among these values, $\vert \visit_{h,h+1} \vert$ can be larger than $2/v$ on at most $vn/2$ occasions.
Thus, (\ref{sanenv inc}) holds. 

\begin{figure}
\centering
\epsfig{file=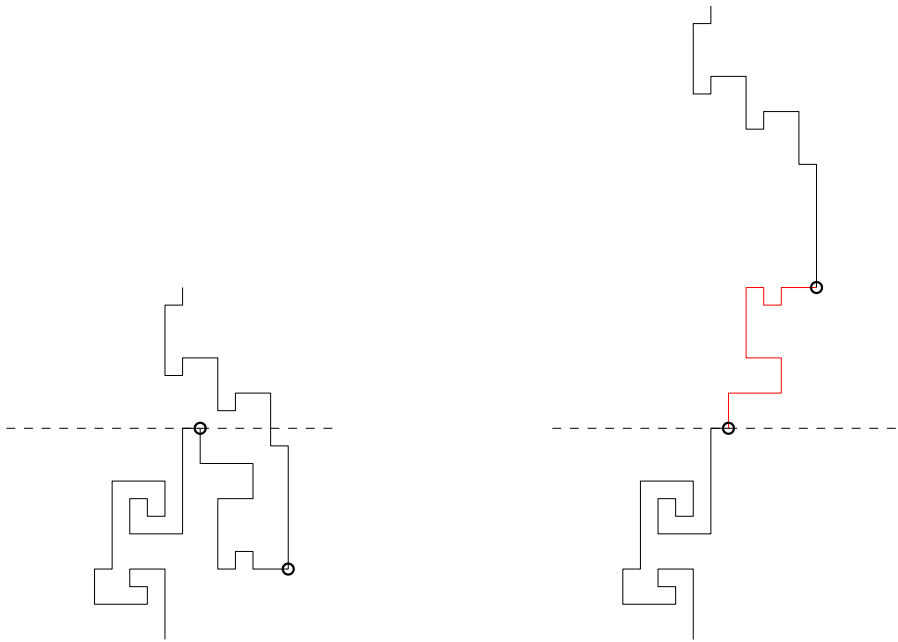, width=12cm}
\caption{In the left-hand sketch, an element of $\sabset_n$ and one of its zigzags is illustrated.  The zigzag's central section is shown in red and the points of zig and of zag are marked by dots. The unfolding of the walk indexed by this zigzag is depicted on the right.}\label{zigzag}
\end{figure}

By (\ref{sanenv inc}), we may successively apply Proposition~\ref{pos dens} with $k$ being set equal to  $\lceil 2/v\rceil$, $\lceil 2/v\rceil-1, \ldots, 3$ and $2$ to conclude that, for any given $v > 0$, there exists some $\delta > 0$, some subsequence $\big\{ v_n: n \in \N \big\}$ of $\big\{ u_n: n \in \N \big\}$ and some  sequence $\big\{ \ep_n: n \in \N \big\}$ with $\ep_n \searrow 0$ such that 
$\vert \sabset_{v_n,v,\delta}^1 \vert \geq e^{-\ep_n n}|\sabset_{v_n,v}|$. 
However, if $\visit_{h,h+1}$ has only one element $e$, then the endpoint $u$ of $e$ with $y(u)=h$ belongs to $\ner_\gamma$, so that $\sabset_{v_n,v,\delta}^1 \subseteq \big\{ \gamma \in  \sabset_{v_n,v}: \vert \ner_\gamma \vert \geq \delta v_n \big\}$, whence follows the theorem. \qed
\bigbreak
We now introduce a concept which plays a central role in the proof of Proposition~\ref{pos dens}. Figure~\ref{zigzag} provides an illustration.

Let $\gamma\in \sabset_n$. A {\em zigzag} of $\gamma$ is a pair $(i,j)$ of indices belonging to $\{0,\ldots,n\}$ that satisfy $i \leq j$ and for which  
\begin{itemize}
 \item the largest of the values $k$ at which the maximum of $y(\gamma_k)$ over $1 \leq k \leq j$ is achieved equals $i$, and
 \item the largest of the values $k$ at which the minimum of $y(\gamma_k)$ over $i \leq k \leq n$ is achieved equals $j$.
\end{itemize}
Let $\zigzag = \zigzag(\gamma)$ denote the set of zigzags of $\gamma$. For $(i,j) \in \zigzag$, we call the image of the subwalk $\gamma_{[i,j]}$ the {\em central section} of $(i,j)$; by the length of this central section, we mean the length of $\gamma_{[i,j]}$, which is $j - i$.
The indices $i$ and $j$ are respectively called the {\em point of zig} and the {\em point of zag} of $(i,j)$. 
%Let $\zig$ and $\zag$ denote the subsets of $\{0,\ldots,n\}$ given by the set of indices which are the point of zig, or the point of zag, of some element of $\zigzag$.  

The following basic properties of the zigzags of a bridge are readily verified. 
 \begin{lemma}\label{lemzigzag}
Let $\gamma \in \sabset_n$. 
\begin{enumerate}
 \item 
 The central sections of the zigzags of $\gamma$ are pairwise disjoint.
 \item If $(i,i) \in \zigzag$, then  $\gamma_i \in \ner_\gamma$.  
\end{enumerate}
 \end{lemma}
We now describe a one-step operation which acts on a bridge and unfolds a given zigzag. To do so, for $v \in \R^d$, we let $\reflect_v :\R^d \to \R^d$ denote the orthogonal reflection with respect to the hyperplane $\big\{ u \in \R^d: y(u) = y(v) \big\}$.
\begin{definition}
For $\gamma \in \sabset_n$ and $(i,j) \in \zigzag$, let  $\unfold{i}{j}{\gamma}  \in \sabset_n$ be given by $\unfold{i}{j}{\gamma} = \gamma_{[0,i]} \circ \reflect_{\gamma_i}\big(\gamma_{[i,j]}\big) \circ \gamma_{[j,n]}$.
\end{definition}

\begin{lemma}\label{lemzzcompose}
Let $\gamma \in \sabset_n$.
\begin{enumerate}
\item Let $(i,j) \in \zigzag$. Then  $\gamma_i$ and $\unf_{(i,j)}(\gamma)_{j}$ each belong to $\ner_{\unf_{(i,j)}(\gamma)}$. 
 \item Let $(i,j) \in \zigzag$. Then $\zigzag(\gamma) \setminus \{ (i,j) \}  \subseteq \zigzag\big( \unfold{i}{j}{\gamma}  \big)$. 
 \item Let $z_1,z_2 \in \zigzag$. Then $\unf_{z_2}\unf_{z_1}(\gamma) = \unf_{z_1}\unf_{z_2}(\gamma)$.
\end{enumerate}
\end{lemma}
The proof of Lemma~\ref{lemzzcompose} is straightforward and omitted. Lemma \ref{lemzzcompose}(3) shows that the order of application of the maps $\unf_z$ is immaterial. This permits us to extend the action of $\unf$ so that several zigzags are undone at once.
\begin{definition}
Let $\gamma \in \sabset_n$ and  $Z \subseteq \zigzag$. Set $\unf_Z(\gamma) \in \sabset_n$ equal to the bridge obtained by iteratively applying  to $\gamma$ the maps $\unf_z$ for $z \in Z$. 
\end{definition}
We finish these preliminaries with a trivial lemma.
\begin{lemma}\label{lemunframp}
Let $\gamma \in \sabset_n$ and  $Z \subseteq \zigzag$. Then $\unf_Z(\gamma) \in \sabset_n$ and $y(\unf_Z(\gamma)_n) \geq y(\gamma_n)$. 
\end{lemma}
We are now in a  position to prove Proposition~\ref{pos dens}.
\medbreak
\noindent{\bf Proof of Proposition \ref{pos dens}.}
The proof is presented as a study of three cases. The walks in $\sabset_{v_n,v,\delta}^k$ are divided into the three following types, and the respective cases consider the event that walks of the given type form a positive proportion of $\sabset_{v_n,v,\delta}^k$:
\begin{enumerate}
 \item {\bf Many renewal points}: walks with a positive proportion of renewal points are common, and there is nothing to prove;
 \item {\bf Few zigzags}: walks with few zigzags are common; here, we unfold {\em all} zigzags at subexponential cost, thereby reducing the cardinality of $\visit_{h,h+1}$ for all levels $h$ for which $\vert \visit_{h,h+1} \vert$ is not already one;
 \item {\bf Many zigzags and few renewal points}: if such walks as these are common, we unfold a collection of short zigzags of such a walk, the number of unfolded zigzags being both a tiny proportion of all zigzags but at the same time of far greater size than the existing number of renewal points; the former bound yields a large choice for such joint unfoldings, while the latter provides for efficient reconstruction of this choice given the outcome. Thus, the multi-valued map principle shows that image walks, which contain many  renewal points because such points are generated wherever unfoldings of zigzags were made, are far more numerous than the preimage walks (which are common by assumption in this case). In spirit, this aspect of the argument has similarities to Kesten's proof~\cite{kestenone} of the pattern theorem for self-avoiding walk.
\end{enumerate}

\medbreak 
\noindent{\bf Case 1: many renewal points.} Suppose that there exists $\delta'>0$ and a subsequence  $\big\{ v_n: n \in \N \big\}$ of $\big\{ u_n: n \in \N \big\}$ such that, for each $n \in \N$, 
\begin{equation*}
\Big\vert \Big\{ \gamma \in \sabset_{v_n,v,\delta}^k: \vert \ner_\gamma \vert \geq \delta' v_n \Big\} \Big\vert \geq \tfrac{1}{4} \big\vert \sabset_{v_n,v,\delta}^k \big\vert   \, .
\end{equation*}
If $\gamma \in \sawset_n$ and $\gamma_i \in \ner_\gamma$ for some $0 \leq i \leq n-1$, note that $\vert \visit_{h,h+1} \vert = 1$ for $h = y(\gamma_i)$. 
Thus if $\gamma \in \sabset_{v_n,v,\delta}^k$ satisfies $\vert \ner_\gamma \vert \geq \delta' v_n$ then $\gamma \in \sabset_{v_n,v,\delta'}^{k-1}$ because $k \geq 2$.
\medbreak
\noindent{\bf Case 2: few zigzags.} Suppose that there exists a sequence $\big\{ \ep_n: n \in \N \big\}$ with $\ep_n \searrow 0$ for which
\begin{equation}\label{eqcasetwo}
\Big\vert \Big\{ \gamma \in \sabset_{u_n,v,\delta}^k: \vert \zigzag_\gamma \vert \le \ep_n u_n \Big\} \Big\vert \geq \tfrac{1}{4} \big\vert \sabset_{u_n,v,\delta}^k \big\vert \, .
\end{equation}
Let $\gamma \in \sabset_n$ and $(i,j) \in \zigzag_\gamma$. If $(i,j)$ is a zigzag satisfying $i < j$ and if $h \in \Z$ satisfies $y(\gamma_j) \leq h< y(\gamma_i)$, then $\vert \visit_{h,h+1}(\gamma') \vert \leq \vert \visit_{h,h+1}(\gamma) \vert - 2$ where $\gamma' = \unf_{(i,j)}(\gamma)$.

If $h \in \Z$ is such that $\vert \visit_{h,h+1}(\gamma) \vert  \geq 2$, then the successive unfoldings performed by $\unf_{\zigzag_\gamma}$ will leave the set $\visit_{h,h+1}$ unchanged except for changing the value of $h$ until the zigzag whose central section includes an edge of $\visit_{h,h+1}(\gamma)$ is unfolded, at which time, the value of $|\visit_{h,h+1}|$ (with the value of $h$ updated) will drop by at least two; subsequently, this set will remain unchanged except for further changes to the value of $h$. Writing $\gamma' = \unf_{\zigzag_\gamma}(\gamma)$ for the outcome of unfolding all of $\gamma$'s zigzags, 
we thus see that the number of $h \in \Z$ such that $1 \leq  \vert \visit_{h,h+1}(\gamma') \vert \leq k-1$ is at least the number of $h \in \Z$ such that $1 \leq  \vert \visit_{h,h+1}(\gamma) \vert \leq k$. 
Also applying Lemma~\ref{lemunframp}, we find that, for $k \geq 2$,
\begin{equation}\label{equnfinc}
\unf_{\zigzag_\gamma} \big( \sabset_{u_n,v,\delta}^k  \big) \subseteq \sabset_{u_n,v,\delta}^{k-1} \, .
\end{equation}
For any $\phi \in \sabset_{u_n}$, we claim that  
\begin{equation}\label{equnfclaim}
\Big| \unf_{\zigzag_\gamma}^{-1}(\phi) \cap \Big\{ \gamma \in \sabset_{u_n}: \vert \zigzag_\gamma \vert \leq \ep_n u_n  \Big\} \Big| \leq  {{u_n} \choose {2\ep_n u_n}} \, . 
\end{equation}
Indeed, if $\phi = \unf_{\zigzag_\gamma}(\gamma)$ for given $\phi$ and unknown $\gamma$, then, to determine $\gamma$, it is enough to know the union of the points in the output that correspond to all of the 
unfolded points of zig and of zag generated as $\unf_{\zigzag_\gamma}$ is formed. This set of points has cardinality $2\ep_n u_n$, whence~(\ref{equnfclaim}).

Applying the multi-valued map principle in light of (\ref{equnfinc}) and (\ref{equnfclaim}), and then using (\ref{eqcasetwo}), we find that
\begin{eqnarray}
 |\sabset_{u_n,v,\delta}^{k-1}| & \ge & 
{{u_n} \choose {2\ep_n u_n}}^{-1}\Big|\Big\{ \gamma \in \sabset_{u_n,v,\delta}^k: \vert \zigzag_\gamma \vert \le \ep_n u_n \Big\}\Big| \nonumber \\
 & \ge & \frac14{{u_n} \choose {2\ep_n u_n}}^{-1} \big|\sabset_{u_n,v,\delta}^{k}\big| \nonumber \\
 & \ge &  \tfrac14 \exp \big\{ - 2 \ep_n u_n \log \big( \tfrac{e}{2 \ep_n} \big) + 1 \big\} \big|\sabset_{u_n,v,\delta}^{k}\big| \, , \nonumber
\end{eqnarray}
where we also used ${m \choose \ell} \leq (m/\ell)^\ell e^{\ell - 1}$ for $m \geq \ell \geq 2$.
\medbreak
\noindent{\bf Case 3: many zigzags and few renewal points.} We now treat the remaining case by supposing that
there exist $\delta'>0$ as well as $\big\{ v_n:n\in\N\big\}$ and $\big\{ \ep_n: n \in \N \big\}$ with $\ep_n \searrow 0$ such that, for all $n$ sufficiently high,
\begin{equation}\label{eqmzfr}
\Big\vert \Big\{ \gamma \in \sabset_{v_n,v,\delta}^k: \vert \zigzag_\gamma \vert \geq 2\delta' v_n \, , \, |\ner_\gamma|\le \ep_n v_n \Big\} \Big\vert \geq  \tfrac{1}{2} \big\vert \sabset_{v_n,v,\delta}^k \big\vert \, .
\end{equation}
Set $\mzfr$ equal to the set on the left-hand side of  (\ref{eqmzfr}).

The central sections of a walk's zigzags being disjoint by Lemma~\ref{lemzigzag}(1), an element
$\gamma \in \sabset_{v_n}$ may have at most $\delta' v_n$ zigzags whose central section has length at least $\lceil1/\delta'\rceil$. A zigzag of $\gamma$ is called {\em short} if its central section has length at most $\lceil1/\delta'\rceil$. Writing  $\shortzz_\gamma$ as the set of short zigzags, note then that $\mzfr \subseteq \big\{ |\shortzz_\gamma|\ge \delta' v_n \big\}$.

Let $\delta''>0$ be a parameter whose value will be fixed later (at a value much less than $\delta'$). We now construct a multi-valued unfolding map called $\multiunf$ that will be defined on the set of  $\gamma \in \sabset_{v_n}$ such that $|\shortzz_\gamma|\ge \delta' v_n$. For such $\gamma$, set 
$$
\multiunf(\gamma)=\Big\{ \unf_Z(\gamma): Z\subseteq \shortzz_\gamma \, , \, |Z|=\lfloor \delta'' v_n \rfloor \Big\} \, .
$$
We now aim to apply the multi-valued map principle. First, if $\phi \in \multiunf(\gamma)$ for some $\gamma \in \sabset_n$ such that $| \shortzz_\gamma | \geq \delta'' n$, then $\vert \ner_\phi \vert \geq  \delta'' n$ because, by Lemma~\ref{lemzzcompose}(1), $\ner_\phi$ contains the set of vertices corresponding to the points of zig and of zag generated by the successive unfoldings which comprise the map whose output is $\phi$ (note that the points of zig and zag of a zigzag can be identical). Hence,
\begin{equation}\label{eqoutputmultiunf}
 \text{if $\phi \in \multiunf(\gamma)$ for some $\gamma \in \mzfr$, then   $|\ner_\phi| \geq \delta'' n$} \, .
\end{equation}
In other words, $\multiunf$ is a multivalued map into $\sabset_{v_n,v,\delta''}^1$.

Second, for each fixed $\gamma \in \sabset_n$, the application from $\zigzag_\gamma$ to $\sabset_n$ which maps $Z$ to  $\unf_Z(\gamma)$ is one-to-one. Hence, if $\gamma \in \mzfr$, then
\begin{equation}\label{eqforwardentropy}
|\multiunf(\gamma)|= {|\shortzz_\gamma|   \choose  \delta'' v_n }\ge {\delta' v_n   \choose  \delta'' v_n } \, .
\end{equation}

The only remaining question to address is the size of preimages. Let $\phi \in \multiunf(\gamma)$ for some $\gamma \in \mzfr$.  We claim that
\begin{equation}\label{equnfoldcomp}
| \ner_\phi | \leq \ep_n v_n + 3\lceil1/\delta'\rceil \delta'' v_n \, .
\end{equation}

Indeed, $\vert \ner_\gamma \vert \leq \ep_n v_n$ by assumption. Furthermore, the number of newly created  renewal points is bounded above by $3\lceil1/\delta'\rceil \delta'' v_n$. Indeed, let $Z \subseteq \zigzag_\gamma$ be such that $\phi = \unf_Z(\gamma)$. 
Consider a single unfolding operation $\unf_{(i,j)}$ applied to a walk $\chi$ (so that $(i,j) \in \zigzag_\chi$). 
Note that in the hyperplane given by $y < y(\chi_j)$, $\unf_{(i,j)}(\chi)$ coincides with $\chi$, while the intersection of $\unf_{(i,j)}(\chi)$ with the hyperplane $y > y(\chi_i) + 2\big(y(\chi_i) - y(\chi_j)\big)$ is merely a translate of the intersection with $\chi$ with $y > y(\chi_i)$. As such, 
 renewal points in $\unf_{(i,j)}(\chi)$ either lie in one of these two hyperplanes, in which case, they have counterparts in $\chi$, or they lie in the slab given by the complement of the union of the hyperplanes. This slab has width $3 \big(y(\chi_i) - y(\chi_j)\big)$, which is at most three times the length $j-i$ of the central section of $(i,j)$. We see then that the number of  renewal points for $\unf_{(i,j)}(\chi)$ exceeds this number for $\chi$ by at most $3(j-i)$. From $j- i \leq  \lceil 1/\delta'\rceil$, we verify that the number of newly created  renewal points satisfies the bound we claimed. From $\vert Z \vert = \delta''n$, we obtain~(\ref{equnfoldcomp}). 
 
Suppose given  $\phi = \unf_Z(\gamma)$ for some unknown  $\gamma \in \mzfr$ and for some unknown subset $Z \subseteq \shortzz_\gamma$ such that 
$\vert Z \vert = \delta'' n$. Recalling Lemma~\ref{lemzzcompose}, 
the data $\gamma$ may be reconstructed from $\phi$ provided that the $\vert Z \vert$ pairs of points corresponding to the points of zig and of zag generated by the successive unfoldings of elements of $Z$ are known. 
Also note that given a certain point $i$ assumed to be the point of zig of some $(i,j)\in Z$, then $j-i$ is the length of the central section and therefore $i\le j\le i+\lceil1/\delta'\rceil$. For this reason, the number of possibilities for the point of zag given the point of zig is bounded by $\lceil 1/\delta'\rceil$.

Hence, the number of pairs $(\gamma,Z)$ such that $\unf_Z(\gamma) = \phi$ is at most the product of
the number of subsets of $\ner_\phi$ of size $\delta'' n$ and the quantity  $ \lceil1/\delta'\rceil^{\delta'' n}$,
the first factor counting possible values of the set of points of zig of $Z$ and the second counting such values for the set of points of zag.
In light of \eqref{equnfoldcomp}, this quantity is at most
\begin{equation}\label{eqbackwardentropy}
  {  (\ep_n+3\lceil1/\delta'\rceil\delta'')v_n \choose \delta'' v_n } \lceil1/\delta'\rceil ^{\delta'' v_n}.
\end{equation}
We omit the proof of the following easy fact.
\begin{lemma}\label{little ineq}
Suppose that $n_1 \geq n_2 \geq m$.  Then $\frac{{ n_1 \choose m}}{{ n_2 \choose m}} \geq \big( \tfrac{n_1 - m}{n_2} \big)^m$. 
\end{lemma}
By setting $n_1 = \delta' v_n$, $n_2 =  (\ep_n+3\lceil1/\delta'\rceil\delta'')v_n$ and $m =  \delta'' v_n$, we may apply Lemma~\ref{little ineq} to bound from below the ratio of the expressions in (\ref{eqforwardentropy}) and (\ref{eqbackwardentropy}) 
provided that
we stipulate that $\delta' \geq \ep_n+3\lceil1/\delta'\rceil\delta''$, a condition which may be realized for all $n$ sufficiently high by an appropriate choice of $\delta''$ since $\ep_n \searrow 0$. Note that the other inequality we need,  $\ep_n+3\lceil1/\delta'\rceil\delta'' \geq \delta''$, is automatically satisfied. 
The resulting lower bound on 
the ratio is
\begin{align*}
 \frac{{ \delta' v_n  \choose  \delta'' v_n }}{{   (\ep_n+3\lceil1/\delta'\rceil\delta'')v_n \choose \delta'' v_n }} \lceil1/\delta'\rceil ^{-\delta'' v_n}  \geq &  \; \bigg( \frac{\delta' -\delta''}{(\ep_n+3\lceil1/\delta'\rceil\delta'')\lceil1/\delta'\rceil}\bigg)^{\delta'' v_n} \, .
\end{align*}
By further imposing on the choice of $\delta''>0$ the requirement that, for all $n$ sufficiently high,
$$
\frac{\delta' -\delta''}{(\ep_n + 3\lceil1/\delta'\rceil\delta'')\lceil1/\delta'\rceil }\ge 2 \, ,
$$ 
we can apply the multi-valued map principle when $n$ is large enough to deduce that
\begin{eqnarray}
&  & |\sabset_{v_n,v,\delta''}^{k-1}| \geq | \sabset^1_{v_n,v,\delta''} | = \big\vert \big\{ \gamma \in \sabset_{v_n,v}: \vert \ner_\gamma \vert \geq \delta'' v_n \big\} \big\vert  \nonumber \\
& \ge & 2^{\delta'' v_n} \Big\vert \Big\{ \gamma \in \sabset_{v_n,v,\delta}^k: \vert \zigzag_\gamma \vert \geq 2\delta' v_n\text{ and } |\ner_\gamma|\le \ep_n v_n \Big\} \Big\vert \, , \nonumber
\end{eqnarray}
whose right-hand side is at least $2^{\delta'' v_n-  1} \big\vert \sabset_{v_n,v,\delta}^k \big\vert$ in view of~(\ref{eqmzfr}). This provides the sought inequality in the third case. \qed
\bigbreak
Recall that $\prne$ denotes the law on $\mathrm{iSAB}$ given by 
$\prne(\gamma) = \mu_c^{-\vert \gamma \vert}$ and $\prneinf$ the law on semi-infinite walks $\gamma:\N \to \lattv$ formed by concatenating infinitely many independent samples of $\prne$. We extend the notion of renewal points to infinite bridges in an obvious way. The set of renewal points is still denoted by $\ner_\gamma$.
\medbreak
\begin{lemma}\label{ren sar}
For each $n \in \N$, we have that
$$
|\sabset_n| \, \muc^{-n} = \prneinf \big( n\in \ner_\gamma \big) \, .
$$
The conditional distribution of $\gamma_{[0,n]}$ under $\prneinf \big( \cdot  \big\vert  n \in \ner_\gamma \big)$ equals $\sab_n$.
\end{lemma}
\noindent{\bf Proof.} Note that 
\begin{eqnarray}
|\sabset_n| \, \muc^{-n} & = & \sum_{\gamma \in \sabset: \vert \gamma \vert = n} \muc^{- \vert \gamma \vert} = \sum_{\ell \geq 0} \sum_{\gamma^{[1]},\ldots,\gamma^{[\ell]} \in \isabset: \sum_{i=1}^\ell \vert \gamma^{[i]} \vert = n} \,\prod_{i=1}^\ell \mu_c^{-\vert \gamma^{[i]} \vert} \nonumber \\
  & = & \prneinf \Big( \exists \, \ell \in \N: \sum_{i = 1}^\ell \vert \gamma^{[i]} \vert = n \Big) = \prneinf \big( n\in \ner_\gamma \big) \, .
  \nonumber
\end{eqnarray}
%\vspace{-5mm}
The latter statement of the lemma follows directly. \qed
%{\flushright{\qed}} 
\medbreak
\noindent{\bf Proof of Corollary~\ref{crucial}.}
Assume that $\mathbb E_{\isabset}\big( |\gamma| \big)=\infty$. Fix $\delta>0$. We may select $k \in \N$ such that
$$
\mathbb E_{\isabset}\big( \min\{k,|\gamma|\} \big)\ge 2/\delta\, .
$$
With $\gamma$ denoting a sample of $\prneinf$, we write $\gamma^{[j]}$ for the $j^{\textrm{th}}$ irreducible  bridge entering in the decomposition of $\gamma$ (translated in such a way that it starts from 0). Define
$$
\textsf{X}_j=\textsf{X}_j(\gamma):= \sum_{i =  1}^j \min\{k,|\gamma^{[i]}|\} \, .
$$
We find that
\begin{equation}\label{eq:ner delta}
\prneinf \big( \big\vert \ner_\gamma \cap [0,n]  \big\vert \geq \delta n \big) \le \prneinf(\textsf{X}_{\delta n}\le n) \, .
\end{equation}
By construction, $\textsf{X}_j$ is a sum of $\delta n$ independent and identically distributed bounded random variables whose common mean is at least $2/\delta$. By the exponential Markov inequality, there exists $\varepsilon>0$ such that, for each $n \in \N$, 
\begin{equation}\prneinf(\textsf{X}_{\delta n}\le n)\le e^{-\varepsilon n}\label{eq:exp dec}.\end{equation} 
However, Lemma~\ref{ren sar} gives
$$
\sab_n \big( \vert \ner_\gamma \vert \geq \delta n  \big) = \prneinf \Big( \big\vert \ner_\gamma \cap [0,n]  \big\vert \geq \delta n \, \Big\vert \, \gamma_n \in \ner_\gamma \Big) \, .
$$
Proposition~\ref{number ramps} implies that
$\vert \sabset_n \vert \geq e^{-C \sqrt n} \vert \sawset_n \vert$ and the classical submultiplicativity of the number of walks together with Fekete's lemma leads to $\vert \sawset_n \vert \geq \mu_c^n$. These two facts in unison with Lemma~\ref{ren sar} yield $\prneinf(\gamma_n \in \ner_\gamma)\ge e^{-C\sqrt n}$. Therefore, \eqref{eq:ner delta} and \eqref{eq:exp dec} imply that
$$
\sab_n \big( \vert \ner_\gamma \vert \geq \delta n  \big) 
\leq e^{- \ep n + C \sqrt n} \, 
$$
for any positive $n$. The quantity $\delta > 0$ being arbitrary, this is in contradiction with Theorem~\ref{thm pos dens}. \qed

\section{Incompatibility of the renewal theory for  bridges}\label{sec:ir} 

In this section, we prove Theorem~\ref{no renewal} by making rigorous the argument outlined after the theorem's statement.
In the presentation, indices are rounded to the nearest integer. For the sake of simplicity, we omit the rounding operation in the notation.
\medbreak
We extend the semi-infinite walk measure $\prneinf$ which we introduced after the statement of Theorem~\ref{no renewal} to doubly infinite walks. Let $\Omega$ be the set of bi-infinite walks $\gamma: \Z \to \lattv$ such that $\gamma_0=0$. Let $\prneinfbi$ denote the law on $\Omega$ formed from $\prneinf$ by iteratively inserting infinitely many independent samples of $\prne$ so that the terminal point of the newly inserted renewal block is the starting point of the walk that is presently formed. Let $\mathcal F$ be the $\sigma$-algebra generated by events depending on a finite number of vertices of the walk. 
 \medbreak
 
 We begin by describing a few properties of the measure $\prneinfbi$.  Let $\mathbf{r}$ be the bi-infinite sequence of integers defined by $\mathbf{r}_0=0$ and $\mathbf{r}_{k+1}=\inf\{j>\mathbf{r}_{k}:\gamma_j\in \ner_\gamma\}$ for each $k\in \mathbb Z$. That is, for $k \not= 0$, $\mathbf{r}_k$ is the index of the $\vert k \vert^{\textrm{th}}$ renewal point, counted northwards or southwards away from the origin according to the sign of $k$. Let $\tau:\Omega\rightarrow \Omega$ be the {\em shift} defined by the formula $\tau(\gamma)_i=\gamma_{i+\mathbf{r}_1}-\gamma_{\mathbf{r}_1}$ for every $i\in \mathbb Z$.
 Let $\sigma$ denote reflection in the hyperplane $\{u\in \R^d:y(u) = 0 \}$. 
 \begin{lemma}\label{prop:prop} The measure $\prneinfbi$ has the following properties.
\begin{itemize}
\item[{\rm (1)}] It is invariant under the shift $\tau$. 
\item[{\rm (2)}] The shift $\tau$ is ergodic for $(\Omega,\mathcal F,\prneinfbi)$.

\item[{\rm (3)}] Under $\prneinfbi$, the random variables $(\sigma\gamma_n)_{n\le 0}$ and $(\gamma_n-\gamma_1)_{n\ge 1}$ are independent and identically distributed.
\end{itemize} \end{lemma}
 
\noindent\textbf{Proof.}
 Part (1) is fairly straightforward. Indeed, the law of $\gamma_{[\mathbf{r}_{-n},\mathbf{r}_n]}$ determines, in the high-$n$ limit, the law of $\gamma$ since we work with the $\sigma$-algebra $\mathcal F$. Now, the laws of $\tau(\gamma_{[\mathbf{r}_{-n-1},\mathbf{r}_{n-1}]})$ and $\gamma_{[\mathbf{r}_{-n},\mathbf{r}_n]}$ are the same by construction (the common law is simply the law of $2n$ consecutively concatenated independent irreducible  bridges). The claim follows by letting $n \to \infty$.
 \medbreak
Part (2) is classical. Let $A$ be a shift-invariant event. We aim to prove that $\prneinfbi(A)\in\{0,1\}$. In order to do so, let $A_n$ be an event depending only on vertices $\gamma_{-n},\ldots,\gamma_n$ such that $\prneinfbi(A_n \Delta A) \le \ep$. By extension, $A_n$ depends only on vertices in $\gamma_{\mathbf{r}_{-n}},\ldots,\gamma_{\mathbf{r}_n}$. Invariance under $\tau$ implies that
\begin{align*}
\prneinfbi(A)&=\prneinfbi(A\cap A)=\prneinfbi[A\cap \tau^{-2n}(A)] \, .
\end{align*}
Moreover,
\begin{eqnarray*}
 & & \Big|\prneinfbi\big(A\cap \tau^{-2n}(A)\big) -\prneinfbi \big( A_n\cap \tau^{-2n}(A_n) \big) \Big| \\
 & \leq & \prneinfbi\big(A \Delta A_n \big)  +  \prneinfbi\big( \tau^{-2n}(A) \Delta \tau^{-2n}(A_n) \big)  \leq 2 \ep \, ,
\end{eqnarray*}
the latter inequality again invoking invariance.
Using the independence between the walk before and after $\mathbf{r}_n$, we find that
$$
\big|\prneinfbi(A)-\prneinfbi(A_n)^2 \big| \le 2\ep
$$
which implies that
$$
 \big| \prneinfbi(A)-\prneinfbi(A)^2 \big| \le 4 \ep \, .
$$
The value of $\ep > 0$ being arbitrary, we conclude that $\prneinfbi(A)=\prneinfbi(A)^2$ and thus that $\prneinfbi(A)\in\{0,1\}$.
\medbreak
Part (3) follows from an explicit construction of the law $\prneinfbi$ which emphasises its near symmetry under $\sigma$.
Noting that any  bridge begins by making a move to the north, let $\isabstar$ denote the law on $\sawset$ formed by removing the first edge of the walk and translating
the resulting subwalk one unit southwards. Note that the support of $\isabstar$ includes the trivial walk of length zero. Note also that $\isabstar$ assigns probability $\muc^{-\vert \gamma \vert-1}$ 
 to all $\gamma$ (of positive length) having the property that $\gamma$ visits only $y$-coordinates comprised between the $y$-coordinates of its starting and ending points, and no such $y$-coordinate is visited exactly once; the residue $ \mu_c^{-1}$ of the probability is assigned to the trivial walk. From this description of $\isabstar$, it is clear that this law is invariant under application of $\sigma$ (and the necessary translation).

As such, $\prneinfbi$ may be constructed as follows. First a north-going edge is placed incident to the origin. An element of $\isabstar$ 
is then concatenated to its northerly endpoint. The construction is then iterated. 
To build the walk in the southerly direction, independent elements of $\isabstar$ and north-south edges are alternately joined so that the terminal point of each is the initial point of the presently constructed walk. 
The invariance of $\isabstar$ under $\sigma$ shows that the two random variables in part(3) each have the law of the concatenation of the sequences formed alternately by independent elements of $\isabstar$ and north-south edges. \qed
%For a finite  bridge~$\gamma$, define $f(\gamma)$ by 
%$$f(\gamma)_n= \sigma\left(\gamma_{|\gamma|-n}-\gamma_{|\gamma|}\right)$$
%for every $0\le n\le |\gamma|-1$. Then, under $\prne$, $\gamma$ and $f(\gamma)$ have the same law. By induction, this implies that $(\sigma(\gamma_n))_{[\textbf{ren}_{-k},0]}$ has the same law as $(\gamma_n-\gamma_1)_{[1,\textbf{ren}_k]}$ for every $k>0$. By construction, they are independent. Letting $k \to \infty$, we obtain the result.
 \medbreak
 Recall the notion of diamond point from Section~\ref{secsketch}. 
 As we explained in the sketch in that section, our goal is now to prove that under the assumption that $\mathbb E_{\rm iSAB}(|\gamma|)<\infty$, a positive density of renewal points are in fact diamond points.
Let $\dr_\gamma$ be the set of diamond points of $\gamma$. Extend the definition of $(\mathbf{r}_k)_{k\ge 0}$ to $\prneinf$.

 \begin{proposition}\label{crucial a}
 If $\mathbb E_{\rm iSAB}(|\gamma|)<\infty$, there exists $\delta>0$ such that 
 $$
\prneinf\left(\liminf_{n\rightarrow \infty}\frac {|\dr_\gamma\cap[0,\mathbf{r}_n]|}n\ge\delta \right)=1 \, .
$$

 \end{proposition}

That the length of the  renewal block has finite mean implies that the bridge under $\prneinf$ is tall and narrow, with height growing linearly and width sublinearly. So the probability under $\prneinfbi$ that the positively (or negatively) indexed walk stays within the cone 
$$\big\{u\in \mathbb Z^d:\ (x+y)(u),(y-x)(u)\ge 0\quad\text{or}\quad(x+y)(u),(y-x)(u)\le 0\big\}$$
 is strictly positive. This separation of future and past occurs at a positive density of  renewal points; when it does, the  renewal point is actually a diamond point. A little care is needed as we make rigorous this argument.
 \medbreak
\noindent{\bf Proof of Proposition~\ref{crucial a}.}
Assume $\mathbb E_{\isabset}(|\gamma|)<\infty$. Let us first prove that $\prneinfbi(\gamma_0\in \dr_\gamma)>0$.

The law of an irreducible bridge is invariant under the reflection with respect to the hyperplane $\{u\in\mathbb Z^d:x(u)=0\}$. Therefore, the law of large numbers applied to $(x,y)(\gamma_{\textbf{r}_n})$ implies that, $\prneinf$-almost surely,  
\begin{equation}\label{eqr}
\left(\frac{x(\gamma_{\textbf{r}_n})}n, \frac{y(\gamma_{\textbf{r}_n})}n\right)\longrightarrow (0,\mu)\quad\text{as $n \to \infty$} \, , 
\end{equation}
 where $\mu$ denotes some strictly positive constant. An irreducible  bridge having finite expected size, we deduce that each of $\inf_{k\ge 0}(x+y)(\gamma_k)$ and $\inf_{k\ge 0}(y-x)(\gamma_k)$ is finite $\prneinf$-almost surely. For $\ell \in \N$, write $\rho_\ell$ for $\prneinf$-probability that both of these random variables is at least  $-\ell$, and choose $K \in \N$ so that $\rho_K > 0$.

We now claim that
\begin{equation}\label{eq:rhorho}
 \rho_0   \ge \mu_c^{-2K} \rho_K \, .
\end{equation}

To verify~\eqref{eq:rhorho}, consider an experiment under which the law $\prneinf$ is constructed by concatenating $2K$ independent samples of $\prne$ to the terminal point of which is concatenated an independent sample of $\prneinf$. If each of the $2K$ samples happens to be a walk of length one (which entails that each is a north-south edge) and if the independent copy of $\prneinf$ realizes 
$$ 
\min\Big(\inf_{k \geq 0} (x+y)(\gamma_k)\ ,\ \inf_{k \geq 0} (y-x)(\gamma_k)\Big) \geq -K \, ,
$$ 
then the constructed sample realizes 
$$
\min\Big(\inf_{k \geq 0} (x+y)(\gamma_k)\ ,\  \inf_{k \geq 0} (y-x)(\gamma_k)\Big) \geq 0 \, .
$$
By the two assertions of Lemma~\ref{ren sar}, the probability that the $i^{\textrm{th}}$ sample of $\prne$ is a north-south edge is $\muc^{-1}$. Thus, the experiment behaves as described with probability $\muc^{-2K} \rho_K$, and we obtain~\eqref{eq:rhorho}.

From \eqref{eq:rhorho} and Lemma~\ref{prop:prop}(3), we deduce that
$$
\delta:=\prneinfbi(\gamma_0\in \dr_\gamma)>0 \, .
$$
The shift $\tau$ being ergodic by Lemma~\ref{prop:prop}(2), we obtain
\begin{equation}\label{eq:doubseq}
\prneinfbi\left(\lim_{n\rightarrow \infty}\frac{|\dr_\gamma\cap[0,\textbf{r}_n]|}{n}= \delta\right)=1 \, .
\end{equation}
From this, the statement of the proposition follows immediately. The equality ``$= \delta$'' becomes an inequality ``$\geq \delta$'' because of the possibility that the negatively indexed semi-infinite bridge under $\prneinfbi$ is responsible for destroying diamond points; in fact, it would be simple to argue that \eqref{eq:doubseq} holds for the law~$\prneinf$, but we have no need of this assertion. 
%In light of the ballistic assumption, Corollary~\ref{crucial} implies that there exists $\delta'>0$ such that
%$$
%\prneinfbi\left(\lim_{n\rightarrow \infty}\frac{\textbf{ner}_n}n=\frac1{\delta'}\right)= 1 \, .
%$$
%We deduce that
%$$
%\prneinfbi\left(\lim_{n\rightarrow \infty}\frac {|\dr\cap[0,n]|}n=\delta\delta' \right)=1.
%$$
%Since avoiding $\gamma_{(-\infty,0]}$ is only relevant for the first ramp renewal point, we get
% $$
%\prneinf\left(\lim_{n\rightarrow \infty}\frac {|\dr\cap[0,n]|}n=\delta\delta' \right)=1.
%$$
%and therefore that
%$$
%\prneinfbi\left(\lim_{n\rightarrow \infty}\frac {|\dr\cap[0,n]|}n=\delta\delta' \,\Big| \,\gamma_0\in\dr\right)=1 \, .
%$$
%The law of $(\gamma_n)_{n\ge 0}$ conditioned on the event that $\gamma_0\in\dr$ is exactly $\prneinf$. Indeed, the conditional measure satisfies the property that 
%$$
%\prne^{\otimes\mathbb N}\Big(\gamma_{[0,n]}=\phi\text{ and }\gamma_n\in \dr \, \Big| \,  \gamma_0\in\dr \Big)=\mu_c^{-n}
%$$ 
%for every $\phi\in\sabset_n$. This property characterizes $\prneinf$. Therefore,
%$$
%\prneinf\left(\lim_{n\rightarrow \infty}\frac{|\dr\cap[0,n]|}n\ge\delta\delta'\right)=1 \, .
%$$ 
%The sign $\ge$ comes from the fact that we forget about the curve $\gamma_{(-\infty,0]}$ which could a priori destroy some ramp renewal points. Proposition~\ref{crucial a} follows from another application of the law of large numbers.
\qed

We now modify long bridges using the stickbreak operation in order to obtain Theorem~\ref{no renewal}. 
\medbreak
\noindent{\bf Proof of Theorem~\ref{no renewal}.}  The proof entails some interplay between infinite and finite bridges: in denoting generic infinite or finite bridges, we write $\phi$ and $\gamma$.

For a finite bridge $\gamma$, let 
$$
\textbf{width}(\gamma)=\max\big\{ x(\gamma_i)-x(\gamma_j):\, 0 \le i,j \le |\gamma|\big\} \, . 
$$
The notion of diamond point may naturally be defined also for finite bridges. We enumerate $\gamma$'s successive diamond points as $\big\{ \mathbf{d}_i(\gamma): 1 \le i \le \vert \dr_\gamma \vert \big\}$.

\medbreak
We aim to prove that $\mathbb E_{\textrm{iSAB}}(|\gamma|) =\infty$. We proceed by contradiction.
Set $\mathbb E_{\textrm{iSAB}}(|\gamma|)<\nu<\infty$ and $0<\ep<\delta/20$, where $\delta$ is given by Proposition~\ref{crucial a}. 
\medbreak
Let $\Omega^+$ denote the set of semi-infinite bridges, so that $\phi \in \Omega^+$ whenever $\phi:\N \to \Z^d$ is a self-avoiding walk for which $y(\phi_i) > 0$ for $i > 0$. 
To any $(\phi,n) \in \Omega^+ \times \N$ for which $\mathbf{r}_n(\phi)$ exists, associate the truncated bridge $\phi^{(n)} : = \phi_{[0,\mathbf{r}_n(\phi)]}$. Note then that $\mathbf{r}_n(\phi)$ coincides with 
$\mathbf{r}_n(\phi^{(n)})$; we write $\mathbf{r}_n$ for the common value in what follows. 
Let $\overline{\sabset}_n=\overline{\sabset}_{n}(\ep)$ denote the set of $\phi \in \Omega^+$ such that
\begin{itemize}
\item[(1)] $\mathbf {r}_n \le  \nu  n$,
\item[(2)] $\textbf{width}(\phi^{(n)})\le \ep n$,
\item[(3)] $\mathbf{d}_{\delta n/2}(\phi^{(n)}) \le \mathbf{r}_n$,
\end{itemize}
and note for later use that the data $\phi^{(n)}$ determines whether these properties hold.

These three properties will be important in implementing the stickbreak operation which was depicted in Figure~\ref{stickbreak} as we outlined our approach. We now prove that
\begin{equation}
\prneinf\big(\overline{\sabset}_n\big)\longrightarrow1\text{ as $n \to \infty$} \, .
\label{eq:ar2}
\end{equation}
From $\mathbb E_{\textrm{iSAB}}(|\gamma|) <\nu<\infty$ and $\mathbb E_{\textrm{iSAB}}[x(\mathbf{t}_\gamma)]=0$, where $\mathbf{t}_\gamma$ is the endpoint of the bridge $\gamma$, we find that
$$
\prneinf\Big(\mathbf{r}_n\le \nu n\text{ and }\textbf{width}(\phi^{(n)})\le \ep n\Big)\longrightarrow 1
$$
thanks to the law of large numbers (which holds since the renewal blocks are independent and identically distributed under $\prneinf$). 
Note that, whenever $\phi^{(n)}$ has at least $i$ diamond points, we have that  $\mathbf{d}_i(\phi^{(n)}) \leq \mathbf{d}_i(\phi)$, because the truncation operation $\phi \to \phi^{(n)}$ can only annihilate (and not create) diamond points.
Thus, Proposition~\ref{crucial a} yields that
$$
\prneinf\big(\mathbf{d}_{\delta n/2}(\phi^{(n)}) \le \mathbf{r}_n\big) \longrightarrow 1 \, .
$$
\medbreak
The convergence \eqref{eq:ar2} implies that we will reach a contradiction if we prove that, for each $n \in \N$,
\begin{equation}\label{eq:littlecontra}
\prneinf\Big(\textbf{width}(\phi^{(n)})> \ep n \Big) \ge \left(\frac{\delta}{20\nu\mu_c}\right)^2\, \prneinf\big(\overline{\sabset}_n\big) \, .
\end{equation}
\medbreak
For $\phi \in \Omega^+$ and $\gamma$ a finite bridge, set $\gamma\triangleleft \phi$ if $\phi_{[0,|\gamma|]}=\gamma$ and if $\phi_{| \gamma |}$ is a renewal point of $\phi$. Lemma~\ref{ren sar} implies that
\begin{equation}\label{eqtrileft}
\mu_c^{-|\gamma|}=\prneinf(\phi\in \Omega^+: \gamma\triangleleft \phi) \, .
\end{equation}

 Let $\overline{\sabset}_n^+ = \big\{ \phi^{(n)}: \phi \in \overline{\sabset}_n \big\}$; that is,  $\overline{\sabset}_n^+$ is  the set of finite bridges $\gamma$ having $n$  renewal points and for which $\gamma \triangleleft \phi$ for some infinite bridge $\phi$ belonging to $\overline{\sabset}_n$. The definition of $\overline{\sabset}_n$ has been chosen so that the membership of any $\phi \in \Omega^+$ in this set is determined by the data $\phi^{(n)}$; thus, we may sum~(\ref{eqtrileft}) to learn that
\begin{equation}\label{sarequiv}
\prneinf\big( \overline{\sabset}_n\big)=\sum_{\gamma\in\overline{\sabset}_n^+}\mu_c^{-|\gamma|} \, .
\end{equation}
We now focus our attention on $\overline{\sabset}_n^+$. Consider $\gamma\in \overline{\sabset}_n^+$. In what follows, we adopt the shorthand $\mathbf{r}_n = \mathbf{r}_n(\gamma)$ and $\mathbf{d}_i = \mathbf{d}_i(\gamma)$. 
For each $i\in \big[\frac\delta {10} n,\frac {2\delta} {10} n\big]$ and $j\in\big[\frac {3\delta}{10}n,\frac{4\delta}{10}n\big]$,
let $\stickbreak_{i,j}(\gamma)$ be given by 
$$
\stickbreak_{i,j}(\gamma) = \gamma_{[0,\mathbf{d}_i]} \circ \tau^{e_1} 
\circ 
\rho_{\pi/2,\mathbf{d}_i} \big(\gamma_{[\mathbf{d}_i,\mathbf{d}_{j}]}\big) \circ \tau^{e_1} \circ 
\gamma_{[\mathbf{d}_{j},\mathbf{r}_n]} \, .
$$ 
Here, $\rho_{\pi/2,\mathbf{d}_i}: \R^d \to \R^d$ is the orthogonal symmetry of $\R^d$ given by a
clockwise rotation of angle $\pi/2$ about $\big(x(\gamma_{\mathbf{d}_i}),y(\gamma_{\mathbf{d}_i})\big)$ in each hyperplane $\{z\in \R^d:z_3=x_3,\ldots,z_d=x_d\}$; and $\tau^{e_1}:\{0,1\} \to \Z^d$ is the length-one walk with $\tau^{e_1}(0) = \bf{0}$ and $\tau^{e_1}(1) = e_1$. 
It is  $\stickbreak_{i,j}(\gamma)$ which was illustrated in Figure~\ref{stickbreak}; the figure shows how this walk is self-avoiding, a property which is  a consequence of the definition of a diamond point.

%We also define the multi-valued map
%$$
%\stickbreak(\gamma)=\Big\{\stickbreak_{i,j}(\gamma)~:~ \textstyle\frac\delta {10}n \le i\le \frac{2\delta} {10} n~,~\frac{3\delta}{10}n\le j\le \frac{4\delta}{10}n\Big\} \, .
%$$
%Set 
%$$
%\stickbreak\big(\overline{\sabset}_n^+\big)=\bigcup_{\gamma\in\overline{\sabset}_n^+}\stickbreak(\gamma) \, .
%$$
%
We now check that  $\stickbreak_{i,j}(\gamma)$ is a  bridge for $i\in \big[\frac\delta {10} n,\frac {2\delta} {10} n\big]$ and $j\in\big[\frac {3\delta}{10}n,\frac{4\delta}{10}n\big]$. 
To do so, it suffices to confirm that for $0 \leq k \leq \mathbf{r}_n + 2$,
\begin{equation}\label{eqxstbr}
  0 \le y\big(\stickbreak_{i,j}(\gamma)_k\big) \le y\big(\stickbreak_{i,j}(\gamma)_{\mathbf{r}_n+2}\big) \, .
\end{equation}
% Indeed, the $x$-coordinate does not raise any difficulty since we rotated by $\pi/2$ in the clockwise direction. 
We prove the second inequality only, the first one being proved in exactly the same way. For $\mathbf{d}_{j} + 2 \leq k \leq \mathbf{r}_n + 2$,
$$y\big(\stickbreak_{i,j}(\gamma)_{\mathbf{r}_n+2} \big) \geq y\big(\stickbreak_{i,j}(\gamma)_{k} \big)$$ by construction. For $\mathbf{d}_{i}\le k\le \mathbf{d}_{j}$,
% and 
%$y\big(\stickbreak_i(\gamma)_{\mathbf{d}_{i}} \big) \geq y\big(\stickbreak_i(\gamma)_{j} \big)$ for $0 \leq j \leq \mathbf{d}_{i + 2\ep n + 2}$.
we have
\begin{eqnarray*}
& &  y\big(\stickbreak_{i,j}(\gamma)_{\mathbf{r}_n+2} \big) ~-~ y\big(\stickbreak_{i,j}(\gamma)_{k} \big) \\
&=  &  \Big( y \big(\stickbreak_{i,j}(\gamma)_{\mathbf{r}_n+2} \big) ~-~ y\big(\stickbreak_{i,j}(\gamma)_{\mathbf{d}_{j}+2} \big) \Big) \\
& & \, \,  + \, \Big( y\big(\stickbreak_{i,j}(\gamma)_{\mathbf{d}_{j} + 2} \big)~-~y\big(\stickbreak_{i,j}(\gamma)_{k} \big) \Big) \\
&\ge & \textstyle\frac{\delta}{10} n-\textbf{width}(\gamma_{[\mathbf{d}_i,\mathbf{d}_{j}]})~\ge~ \big(\frac\delta{10}-\ep\big)n~>~\ep n \, .
\end{eqnarray*}
Regarding the first inequality, we used property (3) of $\overline{\sabset}_n$ and $j\le \frac{4\delta}{10} n$ to find that there are at least $\frac{\delta}{10} n$ diamond points of $\gamma$ whose index lies between $\mathbf{d}_{j}$ and $\mathbf{r}_n$, and thus bounded the first term; in the second inequality, the second term was bounded by using property (2). For $0 \le k \le \mathbf{d}_i$, we get 
\begin{eqnarray*}
&  & y\big(\stickbreak_{i,j}(\gamma)_{\mathbf{r}_n+2} \big) ~- ~y\big(\stickbreak_{i,j}(\gamma)_k \big) \\
&\ge & \Big( y\big(\stickbreak_{i,j}(\gamma)_{\mathbf{r}_n+2} \big)~ -~ y\big(\stickbreak_{i,j}(\gamma)_{\mathbf{d}_{j}+2}\big) \Big) \\
& & \, \, + \, \Big( y\big(\stickbreak_{i,j}(\gamma)_{\mathbf{d}_{j} + 2} \big)~-~y\big(\stickbreak_{i,j}(\gamma)_{\mathbf{d}_i} \big) \Big) \\
& & \, \, + \, \Big( y\big(\stickbreak_{i,j}(\gamma)_{\mathbf{d}_i} \big)~-~y\big(\stickbreak_{i,j}(\gamma)_{l} \big) \Big) \\
&\ge &  \textstyle\frac{\delta}{10} n-\textbf{width}(\gamma_{[\mathbf{d}_i,\mathbf{d}_{j}]})~\ge~ \big(\frac{\delta}{10}-\ep\big)n~>~\ep n \, .
\end{eqnarray*}
That the third term after the first inequality is positive was used.
\medbreak
Let 
$$
\Phi~=~\big[\textstyle\frac{\delta}{10} n,\frac{2\delta}{10}n\big]~\times~ \big[\frac{3\delta}{10}n,\frac{4\delta}{10}n\big]~\times~ \overline{\sabset}_n^+ \, .
$$ We introduce the quantity
$$
S:=\sum_{(i,j,\gamma)\in \Phi}\mu_c^{-|\stickbreak_{i,j}(\gamma)|} \, .
$$
One can express $S$ in terms of $\prneinf\big( \overline{\sabset}_n\big)$. Indeed,
$\vert \stickbreak_{i,j}(\gamma) \vert = \vert \gamma \vert + 2$, and therefore
\begin{equation}\label{eqsinferone}
S = \sum_{(i,j,\gamma)\in \Phi}\,\mu_c^{-|\gamma|-2} = \left(\frac{\delta n}{10\mu_c}\right)^2\sum_{\gamma\in \overline{\sabset}_n^+}\mu_c^{-|\gamma|}= \left(\frac{\delta n}{10\mu_c}\right)^2\,\prneinf\big( \overline{\sabset}_n\big) \, .
\end{equation}
In the last equality, we used \eqref{sarequiv}. 
However, $S$ can be expressed in another way, which we now present. 
\medbreak

First note that  the width of $\stickbreak_{i,j}(\gamma)$ is larger than $\frac\delta {10}n>\ep n$ (since the height of the rotated piece is larger than its number of diamond points). 
Furthermore, every  renewal point of $\stickbreak_{i,j}(\gamma)$ corresponds to a  renewal point of $\gamma$, except for  renewal points possibly created whose index lies between $\mathbf{d}_i$ and $\mathbf{d}_{j}+2$. There are at most $\textbf{width}(\gamma_{[\mathbf{d}_i,\mathbf{d}_j]})\le\ep n$ of them. 
But at least $\frac\delta {10} n$  renewal points (in fact even diamond points) were destroyed in the rotation of the central part, so that $\stickbreak_{i,j}(\gamma)$ contains at most $n+\ep n-\frac\delta{10} n\le n$ renewal points. All together, we find that if $\stickbreak_{i,j}(\gamma)\triangleleft \phi$ holds for some $\phi \in \Omega^+$, then $r_n(\phi) \geq | \stickbreak_{i,j}(\gamma) |+2$ and thus $\textbf{width}(\phi_{[0,\mathbf{r}_n(\phi)]})> \ep n$. We deduce that, for any $(i,j,\gamma)\in \Phi$,
\begin{eqnarray*}
& & \mu_c^{-|\stickbreak_{i,j}(\gamma)|} \\
& = & \prneinf\big( \phi\in \Omega^+:\stickbreak_{i,j}(\gamma)\triangleleft \phi \big) \\
& = & \prneinf\big( \phi\in\Omega^+:\stickbreak_{i,j}(\gamma)\triangleleft \phi\text{ and }\textbf{width}(\phi_{[0,\mathbf{r}_n(\phi)]})> \ep n \big) \, .
\end{eqnarray*}
Therefore,
\begin{eqnarray} 
S &= & \sum_{(i,j,\gamma)\in \Phi}\prneinf\big(\phi\in \Omega^+:\stickbreak_{i,j}(\gamma)\triangleleft \phi\ \text{and}\ \textbf{width}(\phi_{[0,\mathbf{r}_n(\phi)]})> \ep n\big) \nonumber \\
&=&\mathbb E_{\rm iSAB}^{\otimes\mathbb N}\Big(\big| \big\{ (i,j,\gamma)\in \Phi:\stickbreak_{i,j}(\gamma)\triangleleft\phi \big\} \big| \cdot \mathbf{1}_{\{\mathbf{width}(\phi_{[0,\mathbf{r}_n(\phi)]})>\ep n\}}\Big) \nonumber \\
& \le &  (2\nu n )^2\,\prneinf\big( \textbf{width}(\phi_{[0,\mathbf{r}_n(\phi)]})> \ep n \big).\nonumber\\
& \ &  \label{eqsinfertwo}
\end{eqnarray}

The inequality here made use of the fact that, for any given  $\phi \in \Omega^+$, the number of elements $(i,j,\gamma)$ of $\Phi$
that satisfy $\stickbreak_{i,j}(\gamma)\triangleleft\phi$ is at most $(2 \nu n)^2$. To see this bound, note that the solution of this reconstruction problem is  
 uniquely determined once we know the indices in $\phi$ of the westerly and easterly endpoints of the first and second one-edge walks $\tau^{e_1}$ used in the stickbreak operation. Each of these indices is at most $\mathbf{d}_{\delta n /2}(\gamma) +2 \leq \nu n +2 \leq 2 \nu n$ by virtue of $\gamma \in \overline{\sabset}^+_n$; hence the desired bound.

The two inferences \eqref{eqsinferone} and  \eqref{eqsinfertwo} about $S$ together yield~\eqref{eq:littlecontra}.
 This concludes the proof.\qed

 \begin{remark}
 Conjecturally, the law of the length of an irreducible bridge belongs to the domain of attraction of a stable law. This prediction arises from the law having infinite mean, and the scale invariance of the conjectural scaling limit. We refer to \cite{AD} for a discussion on the distribution of the scaled renewal points in dimension two. \end{remark}

\section{Deriving Corollaries~\ref{cor:1} and \ref{cor:2}}\label{seccoroll}
 
Note first that Corollary~\ref{cor:1} is a trivial consequence of Theorem \ref{sub SAW}.
 \medbreak
 \noindent\textbf{Sketch of proof of Corollary~\ref{cor:2}.} We explain how to reduce the proof to a standard application of the Ornstein-Zernike theory developed for subcritical self-avoiding walk by \cite{Ioffe}.
 Fix $K>0$. Instead of working on $\mathcal O_\delta$ with $\delta \searrow 0$, we expand the domain $\calO$ by a factor of $n$ and work on $\mathbb Z^d$. Without loss of generality, assume that $||a-b||=1$. We write $\an$ and $\bn$ for sites of $n\calO\cap \mathbb Z^d$ closest to $na$ and $nb$. By Proposition~\ref{number ramps},
 \begin{equation}\label{eqabxc}
\mathbb P_{(n\calO,\an,\bn,\mu_c^{-1})}(|\gamma|\le Kn) \le\frac{\sum_{k\le Kn} e^{C\sqrt k}\sab_n\big(||\gamma_k||\ge k/K \big)}{\sum_{\gamma\subseteq n\calO: \, \gamma_0=\an \, ,\, \gamma_{|\gamma|}=\bn} \,\mu_c^{-|\gamma|}} \, .
 \end{equation}
 Theorem~\ref{sub SAW} implies that the numerator of~(\ref{eqabxc}) decays exponentially fast. Therefore, it suffices to prove that the denominator does not decay exponentially fast. 
Fixing $\delta>0$, we wish to show that 
\begin{equation}\label{eqnodecay}
\sum_{\gamma\subseteq n\calO: \, \gamma_0=\an \, ,\, \gamma_{|\gamma|}=\bn} \,\mu_c^{-|\gamma|}\ge \exp \big\{ -\delta n \big\}
\end{equation}
for $n$ large enough. In the spirit of the proof of Theorem $A$ of \cite{Ioffe}, it may be argued that, for $z<\mu_c^{-1}$,
$$ 
\sum_{\gamma\subseteq n\calO: \, \gamma_0=\an \, ,\, \gamma_{|\gamma|}=\bn} \,z^{|\gamma|}\ge \exp\Big\{ -\textbf{length}_z(\Omega,a,b)n(1+o_{n\rightarrow \infty}(1))\Big\} \, ,
$$
where $\textbf{length}_z(\Omega,a,b)$ is the length of the geodesic from $a$ to $b$ in $\Omega$ for the norm associated to the correlation length of self-avoiding walk of parameter $z$ (see (1.3) of \cite{Ioffe} for the definition of this correlation length). That this correlation length diverges as $z\nearrow \mu_c^{-1}$ follows from the divergence of $\sum_{\gamma \in \sawset} z^{\vert \gamma \vert}$ in this limit (a fact which we noted in the proof of Lemma~\ref{Kesten relation}). Thus, the parameter $z$ can be chosen close enough to $\mu_c^{-1}$ that $\textbf{length}_z(\Omega,a,b)<\delta$. We thus obtain~(\ref{eqnodecay}). \qed
 
\bibliographystyle{plain}

\bibliography{bridge}

\end{document}